\documentclass{article}
\usepackage[english]{babel}
\input xy
\xyoption{all}
\usepackage{amssymb,amsmath,amsthm}

\newcommand{\R}{\mathbb{R}}
\newcommand{\lie}[1]{\mathfrak{#1}}     

\newcommand{\C}{\mathbb{C}}

\newcommand{\hook}{\lrcorner\,}

\newcommand{\SU}{\mathrm{SU}}

\newcommand{\Gtwo}{\mathrm{G}_2}

\newcommand{\SL}{\mathrm{SL}}
\newcommand{\dfn}[1]{\emph{#1}}

\DeclareMathOperator{\Hom}{Hom}

\DeclareMathOperator{\ad}{ad}

\theoremstyle{plain}
\newtheorem{proposition}{Proposition}
\newtheorem{theorem}[proposition]{Theorem}
\newtheorem{lemma}[proposition]{Lemma}
\newtheorem{corollary}[proposition]{Corollary}
\theoremstyle{definition}

\theoremstyle{remark}

\renewcommand{\Re}{\mathfrak{Re}\,}

\newcommand{\Span}[1]{\operatorname{Span}\left\{#1\right\}}
\begin{document}
\title{Half-flat nilmanifolds}
\author{Diego Conti}
\maketitle

\begin{abstract}
We introduce a double complex that can be associated to certain  Lie algebras, and show that its cohomology determines an obstruction to the existence of a half-flat $\SU(3)$-structure. We obtain a classification of the $6$-dimensional nilmanifolds carrying an invariant half-flat structure.
\end{abstract}

\vskip5pt\centerline{\small\textbf{MSC classification}: Primary 53C25; Secondary 53C29, 17B30}\vskip15pt

An $\SU(3)$-structure on a manifold of real dimension $6$ consists of a Hermitian structure $(g,J,\omega)$ and a unit $(3,0)$-form $\Psi$; since $\SU(3)$ is the stabilizer of the transitive action of  $\Gtwo$ on $S^6$, it follows that a $\Gtwo$-structure on a $7$-manifold induces an $\SU(3)$-structure on any oriented hypersurface. If the $\Gtwo$-structure is torsion-free, meaning that it corresponds to a holonomy reduction, then the $\SU(3)$-structure is half-flat \cite{ChiossiSalamon}; in terms of the defining forms, this means
\begin{equation}
 \label{eqn:HalfFlat}
d\omega\wedge\omega=0, \quad d\Re\Psi=0.
\end{equation}
Conversely, it follows from a result of Hitchin \cite{Hitchin:StableForms} that every compact, real-analytic half-flat $6$-manifold can be realized as a hypersurface in a manifold with holonomy contained in $\Gtwo$, though this is no longer true if the real-analytic hypothesis is dropped \cite{Bryant:NonEmbedding}. Moreover, the $\Gtwo$-structure can be obtained from the half-flat structure by solving a PDE (which turns into an ODE in the  homogeneous case), so that the construction of half-flat structures is indirectly a means of constructing local metrics with holonomy $\Gtwo$. Half-flat manifolds are also studied in string theory (see e.g. \cite{Gurrieri}).

An effective technique to obtain compact examples of half-flat manifolds consists in considering left-invariant structures on a nilpotent Lie group, following \cite{Salamon:ComplexStructures}. Six-dimensional nilpotent Lie algebras are classified in \cite{Magnin}; moreover, each associated Lie group has a uniform discrete subgroup (see \cite{Malcev}), giving rise to a compact quotient, called a nilmanifold. Thus, an $\SU(3)$-structure on the Lie algebra determines an invariant $\SU(3)$-structure on the associated nilmanifold, and vice versa. The nilmanifolds admitting certain special types of half-flat structures have been classified  in \cite{ChiossiSwann,ChiossiFino,ContiTomassini}, and an analogous classification in five dimensions has been obtained in \cite{ContiSalamon}; however, the problem of determining all the nilmanifolds admitting an invariant half-flat structure has been open for some time. A related $5$-dimensional geometry has been studied  in \cite{DeAndres:HypoContact}, leading to examples of half-flat structures on solvable Lie algebras.

\smallskip
In this paper we give a complete classification of nilpotent Lie algebras that admit a half-flat structure. For $24$ out of the $34$ equivalence classes of algebras, we give an explicit example of a half-flat structure; some of these examples are not new. In order to prove that the remaining algebras do not admit a half-flat structure, we introduce an obstruction which also applies to non-nilpotent algebras. Said obstruction is based on the fact that given an $\SU(3)$-structure $(g,J,\omega,\Psi)$ and any covector $\alpha\in T^*_xM$, then
\begin{equation}
 \label{eqn:inequality}
\alpha\wedge J\alpha\wedge\omega^2\neq 0.
\end{equation}
This inequality is linear in $\omega^2$, though the dependence on $\Psi$ is more complicated. However, for structures on a Lie algebra $\lie{g}$, it sometimes happens that the closedness of $\Re\Psi$ implies the existence of a fixed $J$-invariant two-plane $V_1\subset\lie{g}^*$; together with \eqref{eqn:inequality}, this determines a simple obstruction to the existence of a half-flat structure.

In turn, the existence of an invariant two-plane $V_1$ depends on a construction that is somewhat independent of the problem of half-flat classification, and may have other applications. In fact, suppose that $V_2$ is a complement of $V_1$ in $\lie{g}^*$, and
\[d(V_1)\subset \Lambda^2V_1, \quad d(V_2)\subset \Lambda^2V_1\oplus V_1\wedge V_2.\]
To such a splitting $V_1\oplus V_2$, called a \dfn{coherent splitting} in this paper, one can associate a natural double complex, and its cohomology $H^{p,q}$. We prove that if $H^{0,3}$ and $H^{0,4}$ are zero, no half-flat structure exists; more precisely, we show that $V_1$ must be $J$-invariant, and obtain a contradiction from \eqref{eqn:inequality}.

Having introduced this obstruction for general Lie algebras, we go back to the nilpotent case. It turns out that exactly two nilpotent Lie algebras do not admit a coherent splitting. What sets these two algebras apart is their property that $[\lie{g},\lie{g}]$ is not abelian; in standard terminology, their derived length is greater than two. However, a direct argument using \eqref{eqn:inequality} shows that they do not admit a half-flat structure. We conclude that nilpotent $6$-dimensional Lie algebras can be divided into three classes:
\begin{itemize}
\item $24$ algebras have both a half-flat structure and a coherent splitting; for every such splitting $H^{0,3}$ is nonzero.
\item $8$ algebras have no half-flat structure; they have a coherent splitting with $H^{0,3}$ and $H^{0,4}$ equal to zero.
\item $2$ algebras have neither a half-flat structure nor a coherent splitting.
\end{itemize}
In other words, we find that in the nilpotent case the existence of a coherent splitting and the dimension of $H^{0,3}$ are sufficient to determine whether a half-flat structure exists.
\section{$\SU(3)$-structures on nilmanifolds}
In this section we introduce some basic properties of $\SU(3)$-structures, as well as some notation which will be used in the sequel.

Given a $6$-manifold $M$, an $\SU(3)$-structure on $M$ is a reduction to $\SU(3)$ of the bundle of frames, which is given by a Hermitian structure $(g,J,\omega)$ and a unit form $\Psi=\psi^++ i\psi^-$ in $\Omega^{3,0}$, where the $\psi^\pm$ are real.  At each point $x$ of $M$, one can always find an orthonormal coframe
\[\eta^1,\dotsc, \eta^6\in T^*_xM,\]
such that
\begin{equation}
 \label{eqn:SU3}
\omega=\eta^{12}+\eta^{34}+\eta^{56},\quad  \psi^+ + i\psi^-=(\eta^1+i\eta^2)\wedge (\eta^3+i\eta^4)\wedge (\eta^5+i\eta^6).
\end{equation}
Working with this frame, one can easily verify that  $g$ and $J$ are characterized by
\begin{gather*}
2\,\omega(X,JY)\omega^3=g(X,Y)\omega^3=-3(X\hook\omega)\wedge (Y\hook\psi^+)\wedge\psi^+, \quad X,Y\in T_xM.
\end{gather*}
Thus, the pair $(\omega,\psi^+)$ determines the $\SU(3)$-structure, provided the forms are ``compatible'', in the sense that at each point there exists a coframe $\eta^1,\dotsc, \eta^6$ such that \eqref{eqn:SU3} holds. In fact, $J$ depends only on $\psi^+$ and the orientation; Hitchin \cite{Hitchin:TheGeometryOfThreeForms} gives an explicit formula, but all we will need is the following.
\begin{proposition}
\label{prop:J}
Let $(g,J,\omega,\Psi)$ be an $\SU(3)$-structure on $M$. If $v$ is in $T_xM$, $\alpha$ in $T_x^*M$, and
\[ \alpha\wedge (v\hook\psi^+)\wedge\psi^+=0,\]
then \[\alpha(J v)=-(J\alpha)v=0.\]
\end{proposition}
\begin{proof}
We can choose a frame $\eta_1,\dotsc, \eta_6$ at $x$ such that the dual basis of $1$-forms satisfies \eqref{eqn:SU3}.
Since the stabilizer $\SL(3,\C)$ of $\Psi$ acts transitively on $\R^6\setminus\{0\}$, we can assume $v=\eta_1$, so
\[ \alpha\wedge (\eta_1\hook\psi^+)\wedge\psi^+=\alpha\wedge 2\eta^{13456}\]
which vanishes if and only if $\alpha(\eta_2)=\alpha(J\eta_1)=0$.
\end{proof}
The $\SU(3)$-structures we are interested in arise in the following way. A nilmanifold is the quotient $\Gamma\backslash G$ of a nilpotent Lie group $G$ by a uniform discrete subgroup. We say an $\SU(3)$-structure on $G$ is invariant if the forms $(\omega,\psi^+)$ are left-invariant. Such a structure passes to the quotient, giving rise to an $\SU(3)$\nobreakdash-structure on the nilmanifold $\Gamma\backslash G$. We then say that $(\omega,\psi)$ is an \dfn{invariant} structure on the nilmanifold.

It will be convenient to define an $\SU(3)$-structure  on a six-dimensional Lie algebra $\lie{g}$ as a pair \[(\omega,\psi^+)\in\Lambda^2\lie{g}^*\times\Lambda^3\lie{g}^*,\]
such that \eqref{eqn:SU3} is satisfied for some basis $\eta^1,\dotsc,\eta^6$ of $\lie{g}^*$ and some $3$-form $\psi^-$.

Nilpotent Lie algebras of dimension $6$ are classified in \cite{Magnin}; up to equivalence, they are $34$, including the abelian algebra. They are listed in Tables~\ref{table:HalfFlat} and~\ref{table:NoHalfFlat}. By \cite{Malcev}, each of them has a uniform discrete group, i.e. gives rise to a nilmanifold. Thus, giving a nilmanifold with an invariant $\SU(3)$-structure is equivalent to giving a nilpotent Lie algebra with an $\SU(3)$-structure. The aim of this paper is to classify nilpotent Lie algebras with an $\SU(3)$-structure satisfying \eqref{eqn:HalfFlat}.

\section{Coherent splittings}
In this section we introduce a double complex that can be associated to certain $6$-dimensional Lie algebras. We prove that the cohomology groups of this double complex give an obstruction to the existence of a half-flat structure on the Lie algebra. In Section~\ref{sec:classification} we shall use this obstruction to prove that eight nilpotent Lie algebras do not admit a half-flat structure.

Let $\lie{g}$ be a Lie algebra, and suppose $\lie{g}^*=V_1\oplus V_2$ as vector spaces, where $\dim V_1=2$. The algebra of exterior forms on $\lie{g}$ inherits the structure of a bigraded vector space, according to
\[\Lambda^{p,q}=\Lambda^p V_1\otimes \Lambda^qV_2.\]
We shall say $V_1\oplus V_2$ is a \dfn{coherent splitting} if
\begin{equation}
 \label{eqn:coherent}
d(\Lambda^{p,q})\subset\Lambda^{p+1,q}+\Lambda^{p+2,q-1}.
\end{equation}
 One can give the same definition for $V_1$ (and even $\lie{g}$) of arbitrary dimension, but we will only need to consider the case that $\dim V_1=2$.

The relevance of coherent splittings is that they enable us to decompose $d$ into the sum of two differential operators
\[\delta_1\colon\Lambda^{p,q}\to\Lambda^{p+1,q}, \quad \delta_2\colon\Lambda^{p,q}\to\Lambda^{p+2,q-1}.\]
Then
\[0=d^2=(\delta_1+\delta_2)^2\]
implies that $(\Lambda^{p,q},\delta_1,\delta_2)$ is a double complex whose non-zero terms have the following pattern:
\begin{equation}
 \label{eqn:pattern}
\begin{gathered}
 \xymatrix{
\dots\\
\Lambda^{0,3}\ar[u]^{\delta_1}\ar[r]^{\delta_2}&\Lambda^{2,2}\\
& \Lambda^{1,2}\ar[u]^{\delta_1}\\
& \Lambda^{0,2}\ar[u]^{\delta_1}\ar[r]^{\delta_2}& \Lambda^{2,1}\\
&	        & \dots\ar[u]^{\delta_1}\\
}
\end{gathered}
\end{equation}
By construction, the associated total complex is $(\Lambda^*\lie{g}, d)$; in particular each diagonal $\bigoplus_{p+q=k}\Lambda^{p,q}$ has only a finite number of non-zero summands.
In principle one might relabel the indices, by defining
\[A^{h,k}=\Lambda^{h+2k,-k},\]
so that $\delta_1$ becomes of bidegree $(1,0)$ and $\delta_2$ of bidegree $(0,1)$. For greater clarity, we shall use  the ``natural'' indices $(p,q)$ instead.

We now recall a standard construction for double complexes. For each $k$, we have a filtration
\begin{equation}
\label{eqn:Filtration}
\Lambda^{2,k-2}\subset \Lambda^{2,k-2}+\Lambda^{1,k-1}\subset \Lambda^{2,k-2}+\Lambda^{1,k-1}+\Lambda^{0,k}=\Lambda^{k}\lie{g}^*.
\end{equation}
Notice that whilst the spaces $\Lambda^{p,q}$ depends on both $V_1$ and $V_2$, this filtration depends only on $V_1$.
We denote by $Z^k$ the space of closed invariant forms in $\Lambda^k\lie{g}^*$. Taking the intersection with $Z^{k}$,  the  filtration \eqref{eqn:Filtration} determines a filtration
\[Z^{k}_2\subset Z^{k}_1\subset Z^{k}_0=Z^{k},\]
and taking the quotient by the $d$-exact forms, we obtain yet another filtration
\[H^{k}_2\subset H^{k}_1\subset H^{k}_0=H^{k}.\]
We can now define the cohomology groups
\[H^{p,q}=\frac{H^{p+q}_p}{H^{p+q}_{p+1}}.\]
By construction, these groups do not depend on $V_2$ but only on $V_1$, and
\begin{equation}
\label{eqn:GradedCohomology}
H^k=\bigoplus_{p+q=k}H^{p,q},
\end{equation}
where $H^k$ is the $k$-th cohomology group of $(\Lambda^*\lie{g},d)$. Taking dimensions, we can rewrite \eqref{eqn:GradedCohomology} as
\[b_k=\sum_{p+q=k}h^{p,q},\]
with obvious notation.

We are now ready to prove the main result of this section.
\begin{theorem}
\label{thm:obstruction}
Let $\lie{g}$ be a $6$-dimensional Lie algebra with a coherent splitting such that $h^{0,4}$ and $h^{0,3}$ vanish.
Then $\lie{g}$ has no half-flat structure.
\end{theorem}
\begin{proof}
Suppose that $(\omega,\psi^+)$ is a half-flat structure on $\lie{g}$. Then $\psi^+$ is in $Z^3$ and $\omega^2$ is in $Z^4$; moreover, if $J$ is the complex structure induced by $\psi^+$ and the orientation, for every $\alpha$ in $\lie{g}^*$  we have
\begin{equation}
 \label{eqn:NonDegenerate}
\alpha\wedge J\alpha\wedge\omega^2\neq 0.
\end{equation}
One way to see this is to use a coframe $\eta^1,\dotsc, \eta^6$, such that \eqref{eqn:SU3} hold, and verify the statement for $\alpha=\eta^1$; since $\SU(3)$ acts transitively on the sphere in $\R^6$, the claim follows.

We now prove that \eqref{eqn:NonDegenerate} does not hold for any $\alpha\in V_1$. Since the space of exact forms is contained in $\bigoplus_{p\geq 1}\Lambda^{p,q}$, the hypothesis can be rewritten as
\begin{equation}
\label{eqn:Z3Z4}
Z^3=Z^3_0=Z^3_1, \quad Z^4=Z^{4}_0=Z^4_1.
\end{equation}
 It follows that
the bilinear maps
\[Z^3\times Z^3 \to \Lambda^6, \quad (\psi,\phi)\to \alpha\wedge (v\hook\psi) \wedge \phi\]
are identically zero for $\alpha$ in $V_1$ and $v$ in its annihilator $(V_1)^o$. Indeed, $v\hook \psi$ has type $(1,1)+(2,0)$ and $\alpha\wedge\phi$ has type $(2,2)$, so their product is zero.
Hence, by Proposition~\ref{prop:J} we have $J\alpha\in V_1$, and therefore
\[\alpha\wedge J\alpha\wedge\omega^2\in \Lambda^{2,0}\wedge Z^4,\]
which by \eqref{eqn:Z3Z4} is the trivial vector space.
\end{proof}
As a simple example, consider the solvable Lie algebra
\[\lie{g}=(0,12,13,14,15,16).\]
This notation means that $\lie{g}$ is characterized by the existence of a fixed basis $e^1,\dotsc, e^6$ of $\lie{g}^*$ verifying
\[de^1=0, \quad de^i=e^1\wedge e^i, i=2,\dotsc, 6.\]
If we set \[V_1=\Span{e^1,e^2},\quad V_2=\Span{e^3,\dotsc, e^6},\]  we obtain a coherent splitting such that the  operators $\delta_1\colon\Lambda^{0,q}\to\Lambda^{1,q}$ are injective; hence, the $h^{0,q}$ are zero, and by Theorem~\ref{thm:obstruction} no half-flat structure exists on this Lie algebra.

In general, verifying the condition of Theorem~\ref{thm:obstruction} is a straightforward but lengthy computation in linear algebra. We now show how the spectral sequence of $\Lambda^{p,q}$ can be used to speed up the calculations.  We shall denote by $(E_k,\partial_k)$ the $k$-th term of the spectral sequence; with our choice of the indices, $\partial_k$ acts on the graded spaces as \[\partial_k\colon E^{p,q}_k\to E_k^{p+k+1,q-k}.\]
Since $p$ and $p+3$ cannot both be in the range $[0,2]$, the $\partial_2$ are zero and the spectral sequence collapses at $E_2$, as is also clear from  the diagram \eqref{eqn:pattern}.

Recall that a Lie algebra $\lie{g}$ of dimension $6$ is called \dfn{unimodular} if $\ad(X)$ is traceless for all $X$ in $\lie{g}$, or equivalently if $b_6=1$.
\begin{proposition}
\label{prop:spectral}
Let $\lie{g}$ be a unimodular $6$-dimensional Lie algebra with a coherent splitting. Then
\[h^{p,q}=h^{2-p,4-q},\quad h^{0,3}= 1-b_1+b_2-\dim E^{0,2}_1,\]
and $h^{0,4}=0$ if and only if $\Lambda^{2,0}$ is generated by an exact $2$-form.
\end{proposition}
\begin{proof}
We shall fix an isomorphism $\Lambda^{2,4}\cong\R$, inducing isomorphisms
\[\Lambda^{p,q}\cong (\Lambda^{2-p,4-q})^*, \quad \alpha\to \cdot\wedge (-1)^{\left[\frac{p+q}2\right]}\alpha.\]
By hypothesis every $5$-form is closed; since $d=\delta_1$ on $\Lambda^{1,4}$, the space $\Lambda^{2,4}$ has neither $\delta_1$-exact nor $\delta_2$-exact forms. It follows that the diagrams
\[\xymatrix{
\Lambda^{p,q}\ar[d]^\cong \ar[r]^{\delta_1} & \Lambda^{p+1,q}\ar[d]^\cong\\
(\Lambda^{2-p,4-q})^*\ar[r]^{\delta_1^*} & (\Lambda^{1-p,4-q})^*
}
\quad
\xymatrix{
\Lambda^{p,q}\ar[d]^\cong \ar[r]^{\delta_2} &  \Lambda^{p+2,q-1}\ar[d]^\cong\\
(\Lambda^{2-p,4-q})^*\ar[r]^{\delta_2^*} & (\Lambda^{-p,5-q)})^*
 }
\]
commute. Thus, we obtain commutative diagrams
\[\xymatrix{E_1^{p,q}\ar[r]^{\partial_1}\ar[d] & E_1^{p+2,q-1}\ar[d]\\
(E_1^{2-p,4-q})^*\ar[r]^{\partial_1^*} & (E_1^{-p,5-q})^*}\]
where the vertical arrows are isomorphisms. This is sufficient to conclude that $H^{p,q}\cong (H^{2-p,4-q})^*$. By a similar argument, or by applying \eqref{eqn:GradedCohomology}, we also see that
\begin{equation}
 \label{eqn:Poincare}
b_k=b_{6-k}.
\end{equation}
This is actually a known consequence of $\lie{g}$ being unimodular (see e.g. \cite{Greub}).

In order to compute $h^{0,3}$, observe that $E^{1,2}_1=E^{1,2}_\infty=H^{1,2}$, since the  coboundary maps
\[\partial_1\colon E_1^{-1,3}\to E_1^{1,2}, \quad \partial_1\colon E^{1,2}_1\to E_1^{3,1}\]
are zero. By the first part of the proof, $E_1^{2,2}$ has the same dimension as $E_1^{0,2}$. On the other hand, by construction we have
\[\dim E_1^{2,2} - \dim E_1^{1,2} + \dim E_1^{0,2} = \dim \Lambda^{2,2}-\dim\Lambda^{1,2}+\dim\Lambda^{0,2}=0.\]
It follows that $h^{1,2}=2\dim E_1^{0,2}$. Also by construction,
\[\sum_{i=0}^6 (-1)^i b_i = \sum_{i=0}^6 (-1)^i\binom{6}{i}=0\;.\]
Using \eqref{eqn:Poincare}, \eqref{eqn:GradedCohomology}  and the first part of the proof, we obtain
\[2(1-b_1+b_2)=b_3=h^{0,3}+h^{2,1}+h^{1,2}=2h^{0,3}+ 2\dim E_1^{0,2}.\]

Finally, by the first part of the proof  $h^{0,4}$ is zero if and only if $h^{2,0}$ is zero. By construction, $d=0$ on $\Lambda^{2,0}$, so this is equivalent to the one-dimensional space $\Lambda^{2,0}$ being generated by an exact form.
\end{proof}

Thus far, we have considered arbitrary Lie algebras. However, if $\lie{g}$ is nilpotent, the existence of a coherent splitting has some extra consequences, which will be useful in Section~\ref{sec:classification}.  We shall denote by $\ker d$ the kernel of \[d\colon\lie{g}^*\to\Lambda^2\lie{g}^*.\]
Nilpotent Lie algebras are solvable, i.e. $\mathcal{D}^n(\lie{g})=\{0\}$ for some $n$,  where $\mathcal{D}(\lie{g})$ is the derived algebra $[\lie{g},\lie{g}]$; the least such $n$ is called the \dfn{derived length} of $\lie{g}$.
\begin{lemma}
\label{lemma:nilpotent}
Let $\lie{g}$ be a nilpotent Lie algebra of dimension $6$, and let $V_1\oplus V_2$ be a coherent splitting of $\lie{g}$. Then
$V_1\subset\ker d$ and the derived length of $\lie{g}$ is $1$ or $2$.
\end{lemma}
\begin{proof}
Let $\alpha$ be in $V_1$; then $d\alpha$ is in $\Lambda^{2,0}$, and so must be zero because $\lie{g}$ is nilpotent. Hence, $d(V_1)=0$. Together with \eqref{eqn:coherent}, this condition implies
\[[\lie{g},\lie{g}]\subset (V_1)^o , \quad \left[(V_1)^o,(V_1)^o\right]=\{0\},\]
and therefore $\mathcal{D}^2(\lie{g})=\{0\}$.
\end{proof}

We conclude this section with a few remarks concerning the invariance of the numbers $h^{p,q}$. To begin with, we show that the $h^{p,q}$ depend on the choice of $V_1$. Indeed, consider the Lie algebra
\[(0, 0, 0, 0, 12, 13).\]
On this algebra, we find three different coherent splittings  whose cohomology satisfies
\begin{gather*}
 h^{0,3} =2,\;  h^{0,4}=0 \text{ if }  e^{12}\in \Lambda^{2,0},\\
 h^{0,3} =2,\;  h^{0,4}=1 \text{ if } e^{14}\in \Lambda^{2,0},\\
 h^{0,3} =3, \; h^{0,4}=1 \text{ if } e^{23}\in \Lambda^{2,0}.
\end{gather*}
Even the condition $h^{0,3}=0$ depends on $V_1$: indeed, for the Lie algebra
\[(0,0,0,12,13,14)\] we find
\begin{gather*}
 h^{0,3} =0= h^{0,4} \text{ if } e^{12}\in \Lambda^{2,0},\quad
 h^{0,3} =2,  h^{0,4}=0 \text{ if } e^{13}\in \Lambda^{2,0}.
\end{gather*}
Thus, the numbers $h^{p,q}$ are not quite invariants of a Lie algebra, since they depend in an essential way on the choice of the $2$-plane $V_1$. However, it makes sense to remove the assumption on the dimensions of $\lie{g}$ and $V_1$ in the definition of a coherent splitting. Indeed every Lie algebra has a canonical splitting, namely
\[\lie{g}^*=\ker d\oplus V_2,\]
where the choice of $V_2$ is irrelevant. One can then define $\Lambda^{p,q}$ as before, obtaining a double complex if \eqref{eqn:coherent} holds --- equivalently, if $[\lie{g},\lie{g}]$ is abelian. The numbers $h^{p,q}$ obtained in this way are canonical invariants of  solvable Lie algebras with derived length $1$ or $2$.

For nilpotent $6$-dimensional Lie algebras, one can see immediately looking at the list (Tables~\ref{table:HalfFlat} and \ref{table:NoHalfFlat}) that there are only two cases in which $[\lie{g},\lie{g}]$ is not abelian. By the argument used in the proof of Proposition~\ref{prop:spectral}, the $h^{p,q}$ satisfy
\[h^{p,q}=h^{b_1-p,6-b_1-q}.\]
However, these invariants are not relevant for the problem at hand: indeed, the two Lie algebras
\[(0, 0, 0, 12, 13, 24), \quad (0,0,0,12,13,14)\]
have the same cohomology, namely
\[h^{0,2}=0,  h^{1,1}=5, h^{2,0}=1, \quad h^{0,3}=0, h^{1,2}=4,\]
but only the first has a half-flat structure.

\section{The classification}

\begin{table}
 \caption{Nilpotent Lie algebras admitting a half-flat structure}
\label{table:HalfFlat}
\[\begin{array}{|l|l|l|}
\hline
\lie{g} &   \text{Adapted half-flat frame} & \Lambda^{2,0}\\
\hline
0, 0, 12, 13, 23, 14&e^{1},e^{5},e^{2},e^{4},e^{3},e^{6} & e^{12}\\
0, 0, 12, 13, 23, 14 + 25&e^{1}-e^{2},e^{4},e^{5},e^{2},e^{6},e^{3}& e^{12}\\
 0, 0, 12, 13,23, 14 - 25, & e^3,e^6,e^4,e^4-e^2,-e^5,e^1+e^5& e^{12}\\
0,0,12,13,14+23,24+15 & -e^5,e^2,e^4,e^1,\sqrt2(e^3-e^5),\frac1{\sqrt2}e^6 & e^{12}\\
\hline
0, 0, 0, 12, 14, 15 + 23&e^{2}+e^{5},e^{2}+e^{5}+e^{6},e^{4},e^{2},e^{3},e^{1}& e^{12},e^{13}\\
0, 0, 0, 12, 14 - 23, 15 + 34&e^{2},e^{4},e^{3},e^{1},e^{6},e^{5}& e^{13}\\
0, 0, 0, 12, 14, 15&e^{1},e^{3},e^{2},e^{5},e^{4},e^{6}& e^{12},e^{13}\\
0, 0, 0, 12, 23, 14 + 35&e^{1},e^{3},e^{4},e^{5},e^{6},e^{2}& e^{13}\\
0, 0, 0, 12, 23, 14 - 35&e^{1},e^{3},e^{2},e^{6},e^{5},e^{4}& e^{13}\\
0, 0, 0, 12, 13, 14 + 35,    & e^2,e^6,-e^3,e^4,e^1,-e^2-e^5& e^{13}\\
0, 0, 0, 12, 13, 14 +23 &e^{2}-e^{6},e^{1}+e^{5},e^{4},e^{1},e^{6},e^{3}& e^{12},e^{13}\\
0, 0, 0, 12, 13, 24&e^{1},e^{6},e^{2},e^{3},e^{4},e^{5}& e^{12},e^{23}\\
0, 0, 0, 12, 13, 23&e^{1},e^{4},e^{2},e^{5},e^{3},e^{6}& e^{12},e^{13},e^{23}\\
0, 0, 0, 12, 14, 15 + 24,   &e^1,e^3,e^2,e^4,e^3+e^5,-e^6& e^{12}\\
0, 0, 0, 12, 14, 15+ 23+ 24 & e^1,e^3,e^2,e^4-e^2,e^3+e^5,-e^6&e^{12} \\
\hline
0, 0, 0, 0, 12, 14 + 25&e^{1}-e^{6},e^{4},e^{5},e^{2},e^{6},e^{3}& e^{12},e^{24}\\
0, 0, 0, 0, 12, 15 + 34&e^{1},e^{3},e^{5},e^{4},e^{6},e^{2}& e^{13},e^{14}\\
0, 0, 0, 0, 13 + 42, 14 + 23&e^{1},e^{2},e^{3},e^{4},e^{5},e^{6}& e^{12},-e^{14}+e^{23},e^{13}+e^{24},e^{34}\\
0, 0, 0, 0, 12, 14 + 23&e^{1},e^{3},e^{2},e^{4},e^{6},e^{5}& e^{12},e^{13},e^{23}-e^{14},e^{24}\\
0, 0, 0, 0, 12, 13&e^{1},e^{4},e^{2},e^{3},e^{5},e^{6}& e^{12},e^{13},e^{14},e^{23}\\
0, 0, 0, 0, 12, 34&e^{1}+e^{3},e^{1},e^{6},e^{5},e^{2},e^{4}& e^{13},e^{14},e^{23},e^{24}\\
\hline
0, 0, 0, 0, 0, 12 + 34&e^{1},e^{2},e^{4},e^{3},e^{5},e^{6}& e^{13},e^{14},e^{23},e^{24},-e^{12}+e^{34}\\
0, 0, 0, 0, 0, 12&e^{1},e^{3},e^{2},e^{4},e^{5},e^{6}& e^{12},e^{13},e^{14},e^{15},e^{23},e^{24},e^{25}\\
\hline
0, 0, 0, 0, 0, 0&e^{1},e^{2},e^{3},e^{4},e^{5},e^{6}& \text {all $2$-forms}\\
\hline
\end{array}\]
\end{table}

\label{sec:classification}
In this section we classify the nilmanifolds carrying an invariant half-flat $\SU(3)$-structure.
\begin{theorem}
\label{thm:classification}
A  $6$-dimensional nilpotent Lie algebra $\lie{g}$ has a half-flat structure if and only if both  these conditions hold:
\begin{itemize}
 \item $\lie{g}$ has a coherent splitting; and
\item for each coherent splitting of $\lie{g}$, $h^{0,3}\neq 0$.
\end{itemize}
In particular, the nilpotent Lie algebras with (without) a half-flat structure are those listed in Table~\ref{table:HalfFlat} (resp.~\ref{table:NoHalfFlat}).
\end{theorem}
The proof consists of three parts, which we state in the form of lemmas.

\begin{lemma}
\label{lemma:h03nonzero}
Let $\lie{g}$ be one of the nilpotent Lie algebras appearing in Table~\ref{table:HalfFlat}. Then $\lie{g}$ has a half-flat structure and a coherent splitting. Moreover  each coherent splitting of $\lie{g}$ has $h^{0,3}\neq 0$.
\end{lemma}
\begin{proof}
A coherent splitting is determined by a generator of $\Lambda^{2,0}$; by construction, such a generator must be simple, namely the product of two $1$\nobreakdash-forms.
Conversely, a simple non-zero $2$-form $\alpha$ determines a coherent splitting if and only if
\begin{equation}
\label{eqn:alpha}
\alpha\in\Lambda^2\ker d;\quad \alpha\wedge de^i=0, \quad i=1,\dotsc, 6.
\end{equation}
The third column in Table~\ref{table:HalfFlat} gives a basis of solutions for this system of equations for each Lie algebra under consideration. The second column contains an example of a half-flat structure for each algebra, given in terms of an adapted frame $\eta^1,\dotsc,\eta^6$, meaning that in terms of this frame the forms $\omega$ and $\psi^\pm$ can be written as \eqref{eqn:SU3}.

Now let $\alpha$ be a generator of $\Lambda^{2,0}$, and suppose that $h^{0,3}=0$; this means that
\begin{equation}
 \label{eqn:alphaZ3}
\alpha\wedge Z^3=0.
\end{equation}
Since a half-flat structure exists, Theorem~\ref{thm:obstruction} implies that $h^{0,4}>0$. By Proposition~\ref{prop:spectral}, this means that $\alpha$ cannot be  exact. We  distinguish among four different cases.

\emph{i}) If $b_1=2$, the only possibility is $\alpha=e^{12}$, which is exact.

\emph{ii}) If $b_1=3$, using \eqref{eqn:alpha} and the fact that $\alpha$ is not exact, we are in one of the following cases:
\begin{align*}
\lie{g}&=(0, 0, 0, 12, 14, 15 + 23), &  \alpha&=\lambda e^{12}+e^{13};\\
\lie{g}&=(0, 0, 0, 12, 14 - 23, 15 + 34), & \alpha&=e^{13};\\
\lie{g}&=(0, 0, 0, 12, 14, 15), &  \alpha&=\lambda e^{12}+e^{13};\\
\lie{g}&=(0, 0, 0, 12, 23, 14 + 35), & \alpha&=e^{13};\\
\lie{g}&=(0, 0, 0, 12, 23, 14 - 35), & \alpha&=e^{13};\\
\lie{g}&=(0, 0, 0, 12, 13, 24), &  \alpha&=\lambda e^{12}+e^{23},
\end{align*}
where $\lambda$ is a constant. We could now apply Proposition~\ref{prop:spectral}; a faster alternative is to verify that in each case there exists a closed $3$-form $\psi$ with $\alpha\wedge\psi\neq 0$; for instance, we may take $\psi=e^{245}$ for the first five algebras and $\psi=e^{145}$ for the last one. This contradicts our assumption \eqref{eqn:alphaZ3}.

\emph{iii}) If $b_1=4$, for every $\beta$ in $\Lambda^2\ker d$, we have the implications
\[\beta\wedge de^5=0\implies d(\beta\wedge e^5)=0 \implies \beta\wedge \alpha=0,\]
where we have used \eqref{eqn:alphaZ3}. Thus $\alpha$ is a multiple of $de^5$, hence exact.

\emph{iv}) Finally, if $b_1\geq 5$  the linear map
\[\Lambda^2\ker d\xrightarrow{f} \Hom(\Lambda^3\ker d,\Lambda^5\ker d), \quad f(\alpha)(\beta)=\alpha\wedge\beta\]
is injective, contradicting \eqref{eqn:alphaZ3}.
\end{proof}

On the side of non-existence, we have the following:
\begin{lemma}
\label{lemma:h03zero}
Of the nilpotent Lie algebras in Table~\ref{table:NoHalfFlat}, all those with $b_2\geq 3$ admit a coherent splitting with
\[e^{12}\in\Lambda^{2,0}, \quad h^{0,3}=0=h^{0,4}.\]
\end{lemma}
\begin{proof}
The fact that the splitting is coherent is obvious. To determine the cohomology we resort to Proposition~\ref{prop:spectral} and compute the dimension of $E^{0,2}_1$; a basis of $E^{0,2}_1$ for each algebra under consideration is given in the fourth column of Table~\ref{table:NoHalfFlat}, and we see that in each case $h^{0,3}=0$. To prove that $h^{0,4}=0$ we use Proposition~\ref{prop:spectral} and the fact that $e^{12}$ is exact in each case.
\end{proof}

There are two algebras left out, which must be discussed separately, namely the two algebras in Table~\ref{table:NoHalfFlat} with $b_2=2$.
\begin{lemma}
\label{lemma:sporadic}
 Let $\lie{g}$ be one of \[(0,0,12,13,14,34+52), \quad (0,0,12,13,14+23,34+52).\] Then $\lie{g}$ has no coherent splitting and no half-flat structure.
\end{lemma}
\begin{proof}
These algebras have derived length $3$, and therefore no coherent splitting, by Lemma~\ref{lemma:nilpotent}. Let $(\omega,\psi^+)$ be a half-flat structure on $\lie{g}$, and let $\sigma=\omega^2$. Then
\[\sigma\wedge e^{12}=d(\sigma\wedge e^3)\]
is an exact $6$-form, hence zero. By the same token, $\sigma\wedge e^{13}=0$.

We shall prove that
\begin{equation}
\label{eqn:Je1}
Z^3\times Z^3 \to \Lambda^6, \quad (\psi,\phi)\to (e_k\hook\psi)\wedge e^1 \wedge \phi
\end{equation}
is zero for $k=4,5,6$. By Proposition~\ref{prop:J}, this condition implies that $Je^1$ lies in $\Span{e^1,e^2,e^3}$. But then
\[e^1\wedge Je^1\wedge\sigma\in \Span{e^{12},e^{13}}\wedge Z_4=\{0\},\]
giving a contradiction by \eqref{eqn:NonDegenerate}.

Thus, it is sufficient to prove \eqref{eqn:Je1}. Observe that in each case the forms
\begin{equation}
 \label{eqn:FiveExactForms}
e^{123},e^{124},e^{125},e^{134},e^{135}
\end{equation}
are exact, and consequently give zero on wedging with $\psi\in Z^3$; in other terms,
\begin{equation*}
\psi\wedge e^1\in\Span{e^{1234},e^{1235},e^{1236}, e^{1245},e^{1345}}.
\end{equation*}
It follows that for each $k=4,5,6$, the form $e_k\hook \psi\wedge e^1$ is a linear combination of the forms \eqref{eqn:FiveExactForms}; since $\phi$ also belongs to $Z^3$, \eqref{eqn:Je1} follows.
\end{proof}

\begin{table}
 \caption{\label{table:NoHalfFlat}
Nilpotent Lie algebras  admitting no half-flat structure}
\[\begin{array}{|l|l|l|l|l|}
\hline
\lie{g} & b_1 &   b_2 & E^{0,2}_1\\
\hline
0,0,12,13,14+23,34+52 &2& 2& \\
0,0,12,13,14,34+52   &2&  2&\\
\hline
0,0,12,13,14,15 &2&   3&e^{34},-e^{36}+e^{45}\\
0,0,12,13,14,23+15 &2& 3&e^{34},e^{45}-e^{36}\\
0,0,0,12,14,24 &3&5 &e^{34},e^{45},e^{46}\\
0,0,0,12,13+42,14+23  &3&5 &e^{34},e^{45}+e^{36},e^{46}-e^{35}\\
0,0,0,12,14,13+42  &3& 5 &e^{34},e^{45},e^{35}+e^{46}\\
0,0,0,12,13+14,24 &3&  5&e^{34},e^{45}+e^{35},e^{46}\\
0,0,0,12,13,14 &3&  6 &e^{34},e^{35},e^{36}+e^{45},e^{46}\\
0,0,0,0,12,15 &4&  7& e^{34},e^{35},e^{45},e^{56}\\
\hline
\end{array}\]
\end{table}
By the classification of $6$-dimensional nilpotent Lie algebras, this ends the proof of Theorem~\ref{thm:classification}. We conclude by giving a different formulation of the results of this section.

First, we observe that the second part of Lemma~\ref{lemma:h03nonzero} could be proven more uniformly as follows. By a long yet straightforward computation one can check that for each algebra in Table~\ref{table:HalfFlat}, the natural bilinear map
$\Lambda^2\lie{g}^*\times Z^3\to \Lambda^5\lie{g}^*$
induces an injective map
\begin{equation}
\label{eqn:theinjectivemap}\Lambda^2\lie{g}^*\to \Hom(Z^3,\Lambda^5\lie{g}^*).
\end{equation}
This immediately implies that $h^{0,3}\neq 0$.

This argument  can also be used to reformulate Theorem~\ref{thm:classification} in a different way, not involving directly the double complex construction. The essential point is that the map \eqref{eqn:theinjectivemap} fails to be injective for the algebras of Lemma~\ref{lemma:h03zero}; another similar condition has to be introduced to account for the algebras of Lemma~\ref{lemma:sporadic}. We denote by $B^2$ the space of exact $2$-forms.
\begin{corollary}
A  six-dimensional nilpotent  Lie algebra has a half-flat structure if and only if both these conditions hold:
\par\textbullet \quad the natural map $\Lambda^2\lie{g}^*\to \Hom(B^2,\Lambda^4\lie{g}^*)$ is not injective; and
\par\textbullet\quad the natural map $\Lambda^2\lie{g}^*\to \Hom(Z^3,\Lambda^5\lie{g}^*)$ is injective.
\end{corollary}
This statement differs from that of Theorem~\ref{thm:classification} in that the first condition is weakened, whereas the second condition is strengthened. However, using the explicit classification of nilpotent Lie algebras of dimension six, one can check that the two conditions together define the same subset of algebras, thus proving the corollary.

\small\noindent Dipartimento di Matematica e Applicazioni, Universit\`a di Milano Bicocca, via Cozzi 53, 20125 Milano, Italy.\\
\texttt{diego.conti@unimib.it}
\end{document}